\documentclass{elsarticle}
\usepackage{amsmath}
\usepackage{amssymb}
\usepackage{url}
\usepackage{xcolor}
\usepackage{subcaption}
\usepackage{overpic}
\usepackage{multirow}
\usepackage{wrapfig}
\usepackage{tcolorbox}

\definecolor{goldenrod}{rgb}{0.85, 0.65, 0.13}
\definecolor{bluegray}{rgb}{0.5, 0.5, 0.85}
\definecolor{gray}{rgb}{0.5, 0.5, 0.5}
\definecolor{orange}{rgb}{1, 0.65, 0}
\definecolor{darkgreen}{rgb}{0.0, 0.5, 0.0}
\definecolor{lightblue}{rgb}{0.5, 0.5, 1.0}

\definecolor{darkbrown}{rgb}{0.59, 0.29, 0.0}
\DeclareMathOperator{\deter}{det}
\DeclareMathOperator{\df}{d}

\newcommand{\chgA}[1]{#1}

\newcommand{\sss}{polyhedral spline}
\newcommand{\sscn}{polyhedral control net}

\newcommand{\figref}[1]{Fig.~\ref{#1}}

\newcommand{\secref}[1]{Section~\ref{#1}}

\newcommand{\ga}{\alpha}   
\newcommand{\gb}{\beta}

\newcommand{\bc}{\mathbf{c}}   
   
\newcommand{\bp}{\mathbf{p}}

\newcommand{\hl}{highlight line}

\newcommand{\R}{\mathbb{R}}

\definecolor{goldenrod}{rgb}{0.85, 0.65, 0.13}
\definecolor{bluegray}{rgb}{0.5, 0.5, 0.85}
\definecolor{gray}{rgb}{0.5, 0.5, 0.5}
\definecolor{orange}{rgb}{1, 0.65, 0}
\definecolor{lightblue}{rgb}{0.5, 0.5, 1.0}

\definecolor{armygreen}{rgb}{0.13, 0.55, 0.13}

\newcommand{\skp}{\hspace{0.05\textwidth}}

\newcommand{\unitdom}{\Box}
\newcommand{\pa}{\ga}
\newcommand{\bsc}[1]{\mathbf{c}_{#1}} 
\newcommand{\Bps}[2]{\phi^{#1}_{#2}} 
\newcommand{\Bsf}[2]{\psi^{#1}_{#2}}
\newcommand{\bbc}[1]{\mathbf{b}_{#1}} 
\newcommand{\BBf}[2]
{b^{#1}_{#2}} 
\newcommand{\srfp}{\mathbf{p}}
\newcommand{\srfq}{\mathbf{q}}
\newcommand{\splinepiece}{spline piece}
\newcommand{\woblender}[1]{}  

\newcommand{\srf}[1]{\mathbf{x}_{#1}}
\newcommand{\idsrf}{\mathbf{x}^{-1}}
\newcommand{\sol}{\mathtt{u}}
\newcommand{\fh}{\sol_h}
\newcommand{\tr}{\mathtt{t}}   
\newcommand{\st}{\mathbf{s}}

\title{Polyhedral Splines for Analysis\\
\small{special issue for Leszek Demkowicz}}
\author{Bhaskar Mishra and J\"org Peters\\
University of Florida}

\begin{document}
\begin{abstract}
Generalizing tensor-product splines to smooth functions whose control nets outline  topological polyhedra,
bi-cubic \sss s form a piecewise polynomial, first-order differentiable space
that associates one function with each vertex.
Admissible \sscn s consist
of grid-,
star-,
$n$-gon-,
polar- 
and three types of T-junction
configurations.
Analogous to tensor-product splines, \sss s
can both  model curved geometry and 
represent higher-order functions on the geometry.
This paper explores the use of
\sss s for
engineering analysis of curved smooth surfaces
by solving elliptic partial differential equations on free-form surfaces without additional meshing.

\end{abstract}


\maketitle

\section{Introduction}
Discontinuous Galerkin (DG) \cite{Cockburn:2000:DGM}
and Petrov--Galerkin (dPG) methods 
\cite{demkowicz2010class,demkowicz2011class} 
do not require differentiability between elements.
This benefits stability, simplifies  
implementation, enables parallel execution, and even allows locally adapting test function spaces on the fly
for high local approximation order.
By contrast, the widely-used tensor-product spline \cite{Boor:1978:PGS} represents differentiable  polynomial function spaces.
Notably, splines have been used to generalize the
 iso-parametric approach
to higher-order finite elements on curved smooth geometry, see
\cite{Braibant:1984:SOD,Shy:1987:SOD,Au:1993:ISF,Schramm:1993:CGD,Cirak:thin:2000,Hughes2005b} --
but only if the geometry can be outlined by a control net in the form of 
a tensor-product grid: at irregularities --
\begin{wrapfigure}{r}{0.3\linewidth}
\vskip-.25cm
\includegraphics[width=0.9\linewidth]{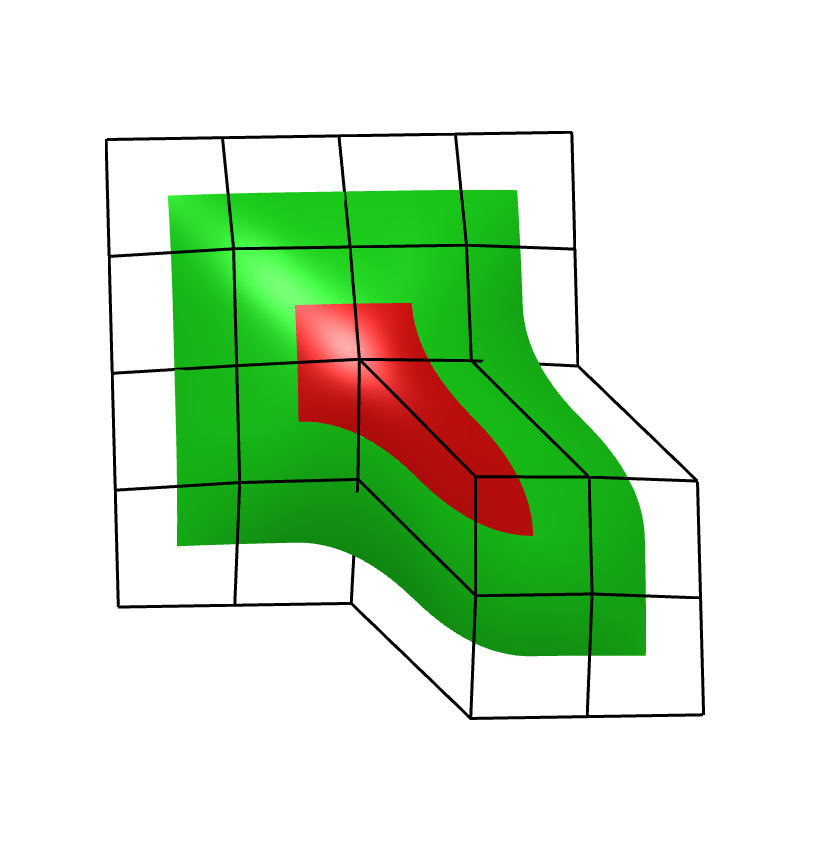}
\vskip-.5cm
\caption{ A \sscn\ and \sss.}
\vskip-.25cm
\label{fig:bdy}
\end{wrapfigure}
 where the tensor-structure breaks down (see the 3-valent and the 5-valent points in \figref{fig:bdy}) -- 
the B-spline and its control net are not well-defined.
Irregularities are no problem for DG spaces.
But when the underlying geometry needs to be both curved and smooth, the DG approach comes up short, 
and requires extra constraints or penalty functions to guarantee a differentiable shape and function
space. Such penalty methods require careful calibration to converge to the proper solution.
Moreover, depending on the differential equation,
a lack of differentiability provides too large a computational space and may so yield outcomes that are not solutions to the original problem,
e.g.\ non-physically discontinuous flow lines.
Differentiability across irregularities is therefore  
both useful and a challenge when devising mathematical software.

Combining differentiability and flexibility,
\sss s \cite{toms:}
extend bi-quadratic (bi-2) tensor-product splines
on regular, grid-like subnets to the non-tensor product sub-net
configurations shown in \figref{fig:lqdnets}.
\def\wid{0.21\linewidth}
\begin{figure}[h]
    \centering
    \begin{subfigure}[t]{\wid}
        \includegraphics[width=\textwidth]{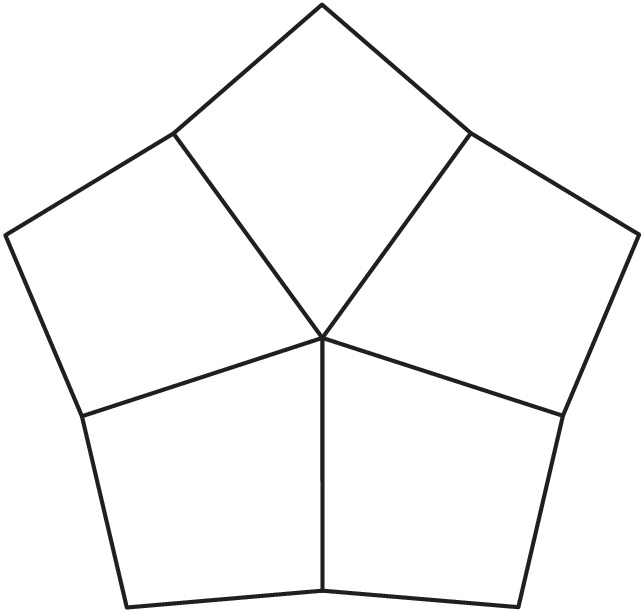}
        \subcaption{star-net $n=5$}
        \label{fig:star-net-n5}
    \end{subfigure} 
    \skp
    \begin{subfigure}[t]{\wid}
        \includegraphics[width=\textwidth]{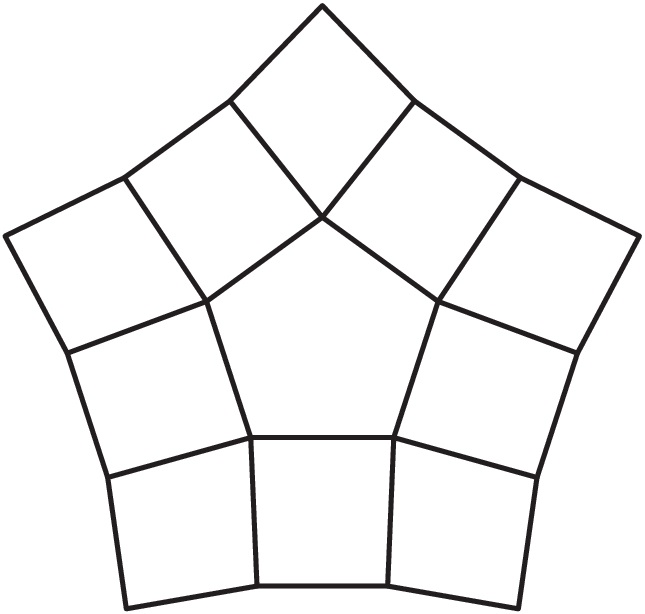}
        \subcaption{$n$-gon-net, $n=5$}
        \label{fig:n-gon-net}
    \end{subfigure}
    \skp
    \begin{subfigure}[t]{\wid}
        \begin{overpic}[scale=.18,tics=10,width=\textwidth]{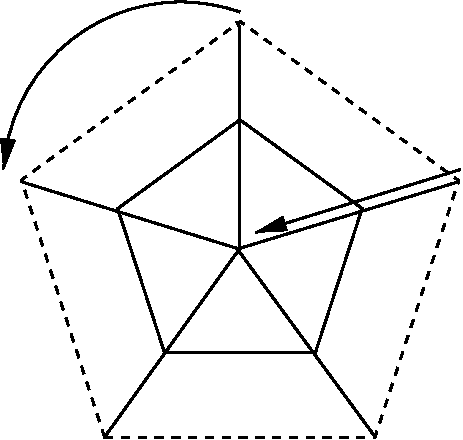}
            \put (10,90) {circular}
            \put (55,55) {radial}
        \end{overpic}
        \subcaption{polar-net $m=5$}
        \label{fig:polar-net}
    \end{subfigure} 
    \bigskip
    
    \begin{subfigure}[t]{\wid}
        \includegraphics[width=\textwidth]{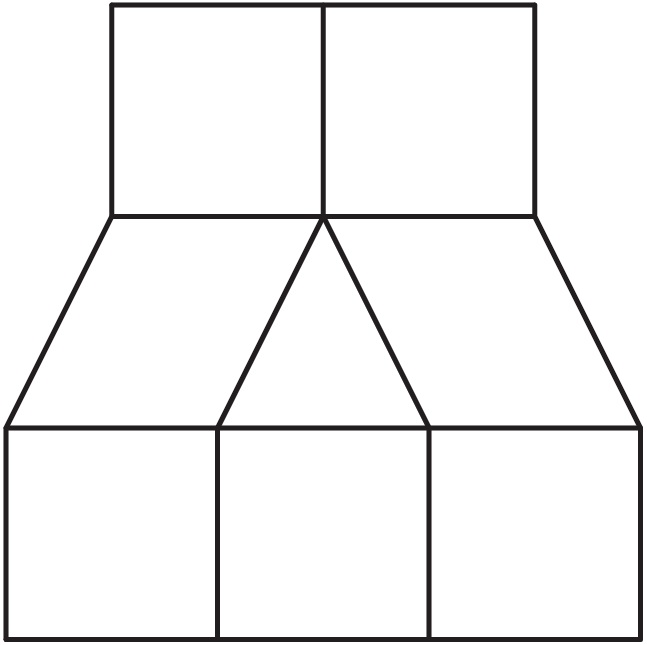}
        \subcaption{$T_0$-net}
        \label{fig:T0-net}
    \end{subfigure} \hspace{0.05\textwidth}
    \begin{subfigure}[t]{\wid}
        \includegraphics[width=1.2\textwidth]{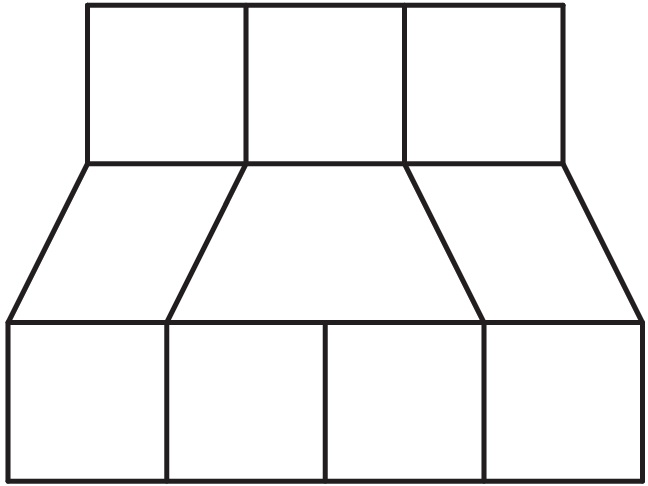}
        \subcaption{$T_1$-net}
        \label{fig:T1-net}
    \end{subfigure}
    \hskip0.1\textwidth
    \begin{subfigure}[t]{\wid}
    \begin{overpic}[scale=.18,tics=10,width=\textwidth]{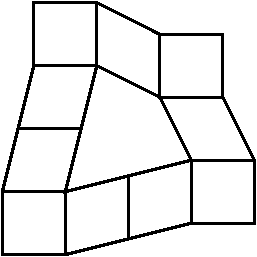}
    \end{overpic}
        \subcaption{$T_2$-net}
        \label{fig:T2-net}
    \end{subfigure}
    \caption{
    Six non-tensor-product \sscn\  patterns. 
    The open source code library \cite{toms:}
    covers $n,m \in \{3,5,6,7,8\}$.
    }
   \label{fig:lqdnets}
\end{figure}
\def\widk{0.15\linewidth}
\def\widm{0.20\linewidth}

\def\wid{0.26\linewidth}
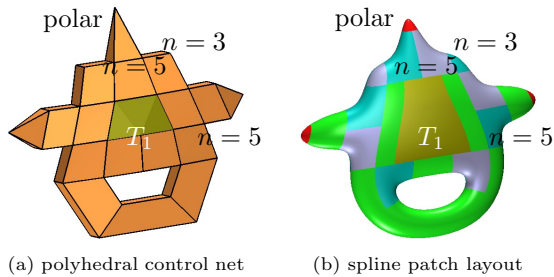
\begin{figure}[h]
\centering
    \begin{subfigure}[t]{\wid}
        \begin{overpic}[scale=.18,tics=10,width=\textwidth]{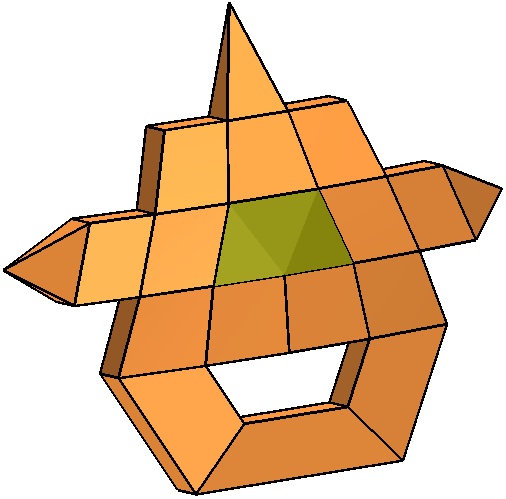}
         \put (15,90) {polar}
         \put (40,70) {\textcolor{black}{${n=5}$}}
         \put (50,40) {\textcolor{white}{$T_1$}}
         \put (65,80) {${n=3}$}
         \put (80,40) {\textcolor{black}{${n=5}$}}
        \end{overpic}
        \subcaption{\sscn\ }
        \label{fig:sscn}
    \end{subfigure}
    \hskip 0.05\textwidth
    \begin{subfigure}[t]{\wid}
        \begin{overpic}[scale=.18,tics=10,width=\textwidth]{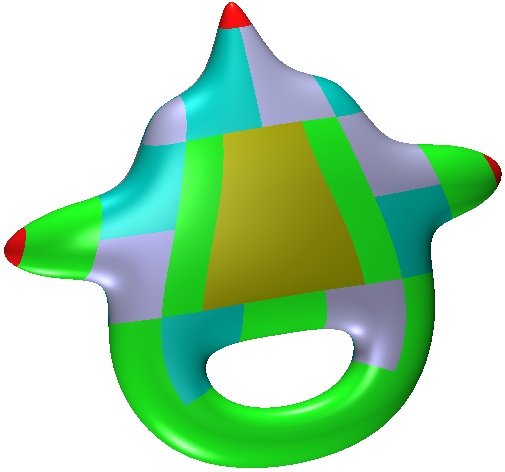}
         \put (15,90) {polar}
         \put (40,70) {\textcolor{black}{${n=5}$}}
         \put (50,40) {\textcolor{white}{$T_1$}}
         \put (65,80) {${n=3}$}
         \put (80,40) {\textcolor{black}{${n=5}$}}
        \end{overpic}
        \subcaption{spline patch layout}
        \label{fig:spline patch layout}
    \end{subfigure}
\caption{`Bottle opener' example of a \sss\  combining non-tensor product patterns in close proximity.
   (a)  \sscn\ with a tight 
   layout: 
   the top pole node is a direct neighbor of four nodes of valence $n=5$.
   5-valent nodes are direct neighbors
   (b) Surface layout including regular
   bi-quadratic \textcolor{green}{bi-2} splines, 
   $n$-sided star-configurations (\textcolor{lightblue}{blue} or \textcolor{gray}{gray}),
   bi-3 \textcolor{red}{polar} caps, and surface pieces covering $T_1$-junctions.
   See also \figref{fig:bottleheat}.
}
\label{fig:tightEnsemble}
\end{figure}

Just like the tensor-product control net, the \sscn\ expresses the neighbor relations of the \sss\  
\chgA{generating} functions, outlines shape and provides handles for manipulating the shape:
the \sscn\ vertices can be used as computational degrees of freedom, 
 for least squares fitting, computing moments or for solving partial differential equations.
The irregular polyhedral patterns listed in \figref{fig:lqdnets} 
can be in close proximity, enabling complex layouts such as \figref{fig:tightEnsemble}.
A \sss\ joins its bi-cubic (bi-3) pieces as a smooth piecewise polynomial function or surface
expressed in Bernstein-B\'ezier form \cite{deboor87e}.
To introduce \emph{creases or discontinuities} 
one can insert facets of zero area, e.g.\ placing opposing sides of a quad onto each other, 
or 
manipulate individual polynomial coefficients of the output.
Due to the sharp degree bound proven in \cite{Karc:2020:LDG}, there exists {no nested} \sscn-based refinement (subdivision algorithm) for bi-3 \sss s,
but the splines 
can be refined by de Casteljau's 
algorithm \cite{farin-curves-and-book-88}
applied per piece, or, non-nestedly, by control net subdivision \cite{catmull78a}.

\medskip
\noindent
\textbf{Overview.}
After brief overview of alternative smooth polynomial function spaces, 
\secref{sec:TNP} reviews tensor-product splines, \sscn s and \sss s.
\secref{sec:compute} extends
the implementation of \sss s \cite{toms:}
to solve
basic partial differential equations
on non-tensor-product meshes -- without additional meshing.

\subsection{Alternative geometric functions spaces}
Commonly used computational polynomial spaces 
for unstructured layout on planar domains include
splines on triangulations \cite{lai_schumaker_2007} 
and radial basis functions \cite{Buhmann:radial}. 
%
For modeling complex geometric free-form shapes, tensor-product spline (NURBS)  domains are carved up into complex regions by a restriction of the domain known as \emph{trimming} (see e.g.\ \cite{Marussig:} in the context of isogeometric design).
Trimming leads to a plethora of downstream
challenges due to heterogeneity in size, parameter orientation, continuity and polynomial degree: algebraic, non-rational pre-image curves can result in gaps in the geometry and the complex domains require special integration rules for engineering analysis.
The animation industry has instead adopted subdivision surfaces \cite{derose:subdiv} that consist in theory of an infinite sequence of nested surface rings,
in practice approximated by a fine faceted model.
The survey of splines for  irregularities on irregular meshes \cite{Peters:2019:SMAI} characterizes
subdivision surfaces as one of three
singular surface constructions:
singularities at corners
\cite{Peters:1991:SIM,Reif1998TURBSTopologicallyUR,Nguyen:2016:RSE,
wu_mourrain_galligo_nkonga_2017},
singular edges
\cite{journals/cad/MylesP11,toshniwal2017multi}
and contracting faces, a.k.a.\
subdivision algorithms \cite{PetersReif:2008:SS}.
Splines leveraging geometric continuity,
i.e.\ differentiable after a local change of variables
include \cite{Peters:1995:SS,Nguyen:2016:CFE,collin2016analysis,blidia:hal-02480959} and for
data fitting and simulation
on planar domains
\cite{Bercovier:2017,journals/cad/KaplST18,journals/cagd/KaplST19,Kapl:smai}.
Analogous to the bi-2 generalizing \sss s in this paper, extending bi-3 tensor-product splines to \sscn s
 and surfaces of good shape
is possible using piecewise polynomials of degree bi-4 or higher,
see e.g.\
\cite{karvciauskas2016generalizing, karvciauskas2019high}.

\section{Tensor-product splines, \sscn s and  \sss s}
\label{sec:TNP}
This section defines and summarizes the evolution of 
splines on regular tensor-product grids,
\secref{sec:tp},
to the \sscn s,
\secref{sec:patterns},
of \sss s,
\secref{sec:sss}.

\subsection{Tensor-Product Splines}
\label{sec:tp}
A tensor-product spline in two variables, $(u,v)$ is a piecewise polynomial function of the form
\cite{Boor:1978:PGS} :
\begin{equation}
    \bp: (u,v) \to
    \bp(u,v) := \sum^{k_1}_{i=0} \sum^{k_2}_{j=0} \bsc{ij} \Bsf{d_1}{i}(u) \Bsf{d_2}{j}(v).
\label{eq:tpspline}
\end{equation}
The coefficients $\bsc{ij}$ scale the B-splines $\Bsf{d_\ell}{i}$ of degree $d_\ell$
and act as degrees of freedom, say for finite element computations or to outline shape. 
Connecting the \emph{control points} $\bsc{ij}$
 to $\bsc{i+1,j}$ and $\bsc{i,j+1}$ wherever possible yields the \emph{control-net}. 
 Due to variation diminishing property and the convex hull property,  the control-net outlines the graph of the spline function, respectively the geometric shape of the spline surface.
Designers edit the control net to shape a surface while automatically maintaining the desired smoothness.

\def\wid{0.25\linewidth}
\def\widb{0.30\linewidth}
\begin{figure}
\centering
    \begin{subfigure}[t]{\widb}
        \includegraphics[width=\textwidth]{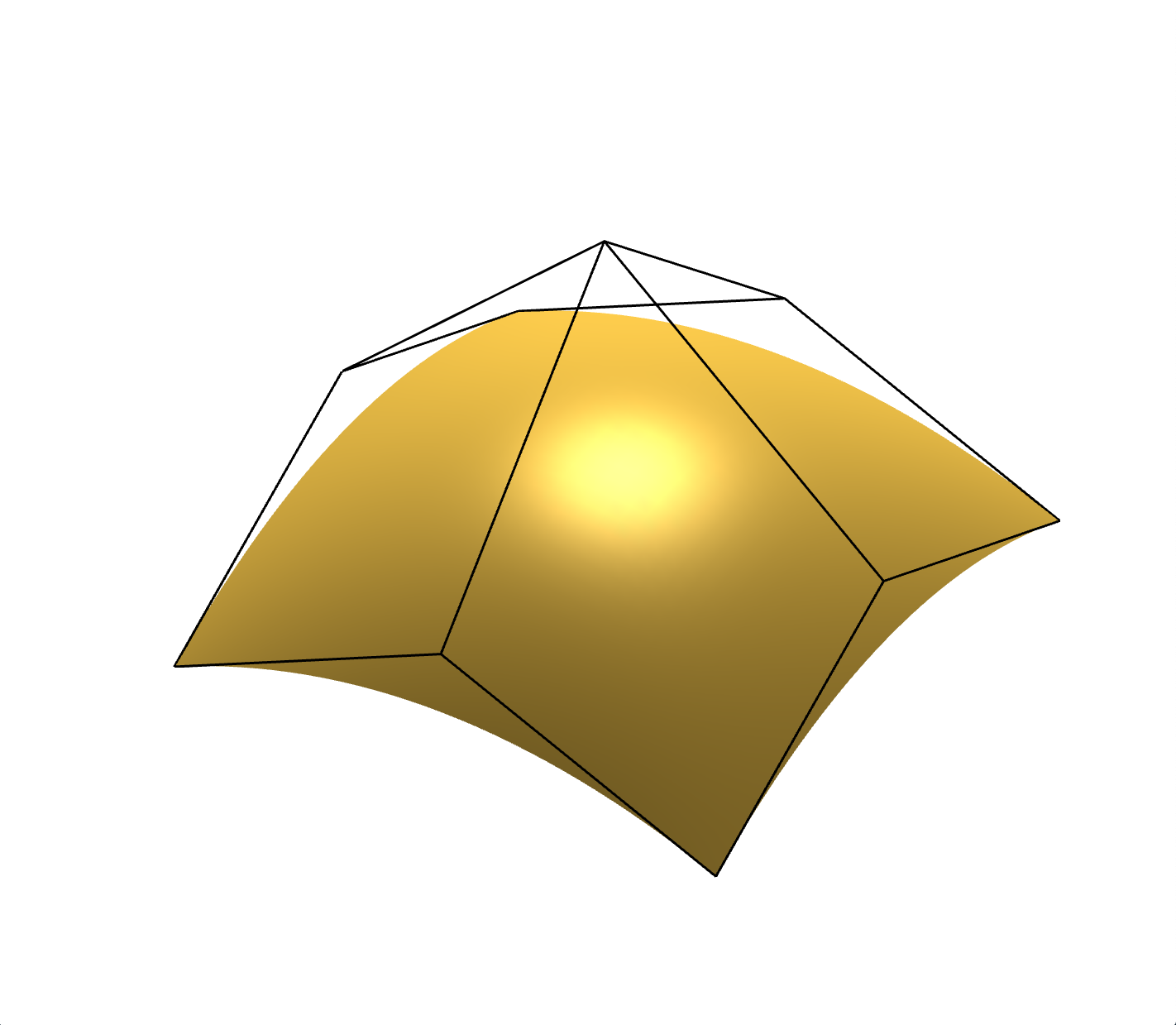}
        \subcaption{BB-net and bi-2 patch}
        \label{fig:control polygon}
    \end{subfigure}
    \begin{subfigure}[t]{\wid}
    \begin{overpic}[scale=.18,tics=10,width=\linewidth]
    {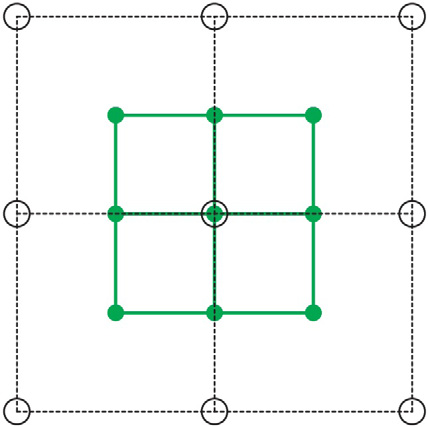}
    \put (116,46) {\textcolor{green}{$\Rightarrow$}}
     \end{overpic}
        \subcaption{bi-2 spline nets}
        \label{fig:bi-2b2bb}
    \end{subfigure}
    \hskip 0.1\textwidth
    \begin{subfigure}[t]{\wid}
        \includegraphics[width=\textwidth]{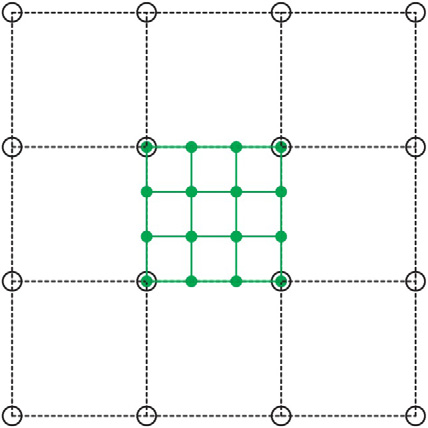}
        \subcaption{bi-3 spline nets}
        \label{fig:bi-3b2bb}
    \end{subfigure}
\caption{
B-spline control points $\circ$ \cite{Boor:1978:PGS} form a control net and BB-coefficients 
\textcolor{darkgreen}{$\bullet$}
\cite{deboor87e} form a BB-net of
$(d+1)^2$ nodes
for a polynomial piece of degree bi-$d$.
\textcolor{green}{$\Rightarrow$} indicates degree-raising for the BB-form.}
\label{fig:B2BB}
\end{figure}

Tensor-product spline control nets form grids of $k_1\times k_2$ quadrilateral faces.
Specifically for our setup, any three by three sub-net can be interpreted as the control 
net of the tensor-product of  B-splines $\Bsf{2}{}(t-k)$ of degree 2 with a \emph{uniform knot sequence}
\cite{deboor87e} 
that define one bi-quadratic (bi-2) polynomial piece:
\begin{equation}
    \bp: (u,v) \to
    \bp(u,v) := \sum^2_{i=0} \sum^2_{j=0} \bsc{ij} \Bsf{2}{}(u-i) \Bsf{2}{}(v-j).
    \label{eq:tpbitwo}
\end{equation}
\subsection{Patterns in \sscn s}
\label{sec:patterns}
\figref{fig:lqdnets} lists irregular  sub-nets of a \sscn\ 
 and Table \ref{tab:sss}  succinctly summarizes
 the structure of the \emph{polyhedral layout configurations} and the resulting \sss s.
Note that for almost all meshes a single Catmull-Clark
\cite{catmull78a} step guarantees a 
quad mesh with star configurations only, i.e.\ a net that controls a \sss\ surface.

Meshes consisting of quadrilaterals are popular in polyhedral 3D modeling
where quad-strips follow the principal directions and delineate features. 
Merging $n$ directions forms either a star-net surrounding an
extraordinary point of valence $n$,
see \figref{fig:lqdnets}a, or an $n$-gon, see 
\figref{fig:lqdnets}b.
Note that overlapping star configurations are admissible and that a quadrilateral face may have multiple non-4-valent vertices.
\figref{fig:tightEnsemble} demonstrates that irregularities can be placed in close 
proximity. 
Often valences $3$ and $5$ suffice for modeling
since vertices or faces with valencies $6$, $7$, $8$ can 
be split and the effect distributed by local re-meshing. However,
it is convenient to also have the freedom to include $n=6,7,8$ to avoid re-meshing
where multiple surface regions meet.
A polar configuration consists of triangles joining at the pole node of high valence to cap off cylindrical 
structure when modeling finger tips and airplane nose cones.
The triangles  are interpreted as quadrilaterals with one edge collapsed 
hence modeled by smoothly joining spline pieces with a (removable) singularity at the pole.

Where two finer quadrilaterals meet a coarser quadrilateral face a T-joint results.
This T gives the 
name to the $T_1$-gon (formally a pentagonal face). A $T_2$-gon (formally a hexagon)
combines two quad-strips at two T-junctions.  
The $T_0$-gon similarly merges neighboring
quad-strips without an explicit T-junction.

\subsection{Piecewise polynomial \sss s}
\label{sec:sss}

\begin{table} 
\caption{Configurations of \sscn s and
\sss s.
All patches are in BB-form of degree bi-3.
A restriction to $n\in \{3, 5,\ldots,8\}$
is only in the distributed code \cite{toms:}; the underlying theory allows for higher $n$.
        }
\begin{center}
\resizebox{\textwidth}{!}{
\begin{tabular}{cccccc}
\hline
\multicolumn{3}{c}{configuration} &  &  &
\\
name & center & surrounded by & Fig. & \# patches & ref
\\
\hline
\\
tensor & $4$-valent point & 4 quads & \ref{fig:B2BB}b & 1 & \cite{deboor87e}
\\
star & $n$-valent  &  $n$ &\ref{fig:star-net-n5} 
&  $n$ for $3,5$ & \cite{KKJP:2014:SMS}
\\
&   &   &
&  $4n$ for $n>5$ & 
\\
$n$-gon &  $n$-gon  & $2n$ & \ref{fig:n-gon-net} 
&  $n$ for $3,5$ &  \cite{KKJP:2014:SMS}
\\
&   &  & 
&  $4n$ for $n>5$ & 
\\ 
$T_0$ &  triangle$^\dag$ & $7$ & \ref{fig:T0-net}
 & $2 \times 2$ &  \cite{sigT:}
\\
$T_1$ &  pentagon$^{\dag\dag}$ & $9$ & \ref{fig:T1-net} & $4 \times 2$ & \cite{sigT:}
\\
$T_2$ &  hexagon$^{\dag\dag\dag}$ & $9$ & \ref{fig:T2-net} & $4 \times 4$ & \cite{sigT:}
\\
polar &  $n$-valent & $n$ triangles & 
\ref{fig:polar-net} & $n$ degenerate & \cite{Karc:2020:SPC}
\\
\hline
\\
\multicolumn{5}{l}{
\scriptsize{
$^\dag$  two vertices of valence 4 and one of valence 5.
$^{\dag\dag}$
four vertices of valence 4 and one of valence 3}}
\\
\multicolumn{5}{l}{
\scriptsize{
$^{\dag\dag\dag}$
three consecutive vertices of valence 4 and 
    two of valence 3
    separated by one vertex of valence 4
         }}
\end{tabular}
}
\end{center}
\label{tab:sss}
\end{table}
A \sss\  is a collection of smoothly-joined 
polynomial pieces in Bernstein-B\'ezier form (BB-form,
\cite{deboor87e,farin-curves-and-book-88})
\begin{equation}
    \srf(u,v) := \sum^{d_1}_{i=0} \sum^{d_2}_{j=0} \bbc{ij} \BBf{d_1}{i}(u) \BBf{d_2}{j}(v), \ \ \ (u, v) \in [0..1]^2
    \label{eq:bbform}
\end{equation}
where $\BBf{d}{k}(t) := \binom{d}{k}(1-t)^{d-k}t^k \in \R$ are the 
Bernstein polynomials of degree $d$ 
and $\bbc{ij}$ are the BB-coefficients. 
For surfaces in 3-space, $\bbc{ij} \in \R^3$
and $\srf{}$ is a piece of the surface, called a \emph{patch}.
For example, a bi-3 patch has $4\times4$ BB-coefficients as in \figref{fig:B2BB}c.
Connecting $\bbc{ij}$ to $\bbc{i+1,j}$ and
$\bbc{i,j+1}$ wherever possible yields the \emph{BB-net}.
The BB-net is usually finer than the B-spline control net of the same polynomial pieces.
(The BB-net represents
a single polynomial piece whereas the B-spline control net represents a \emph{piecewise} polynomial function.
However,
Any B-spline can be expressed as multiple pieces of polynomials in BB-form and any basis function of the BB-form can be expressed in B-spline form with suitably repeated knots \cite{de1986b,deboor87e}.)
The BB-form is evaluated via de Casteljau's algorithm \cite{enwiki:1043959654}, differentiated exactly by forming differences of the BB-coefficients and exactly integrated by forming sums 
\cite{farin-curves-and-book-88}.
To obtain a uniform degree bi-3 in all cases, any patch of degree lower than bi-3 can be expressed as a patch of degree bi-3
by a process called degree-raising \cite{deboor87e,farin-curves-and-book-88}.
Degree bi-3 is therefore the default for the \sss s
used in this paper.


\begin{wrapfigure}{r}{0.4\linewidth}
\vskip -.5cm
\begin{overpic}[scale=.18,tics=10,width=\linewidth]{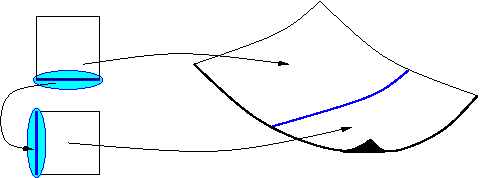}
    \put (40,27) {$\srf{}$}
    \put (40,9) {$\srfq$}
    \put ( 0,10) {$\gb$}
    \put ( 5,-3) {$u$}
    \put (22,19) {$u$}
    \put (65,15) {$E$}
\end{overpic}
\caption{Geometric smoothness between patches $\srf{}$ and $\srfq$ to form an atlas.}
\label{fig:gc}
\end{wrapfigure}
The polyhedral \splinepiece s join by default without gaps or overlap and
with matching first derivatives -- possibly after a reparameterization (a change of variables) as is appropriate for manifolds.
Focusing on the shared boundary $E$ between two
abutting patches $\srfp$ and $\srfq$, see \figref{fig:gc}.  
Let $\gb: \mathbb{R}^2 \to \mathbb{R}$, $ \gb(u) := (u+b(u)v, a(u)v)$
be a (local) reparameterization 
and
$\srfp(u,0) = E = \srfq(u,0)$.
We say that $\srfq$ and $\srfp$ join $G^1$ if the partial derivatives lie in the same plane
and the transversal $v$-derivatives of the two patches lie on opposite sides with respect to
the $u$-derivative along the shared boundary:
\begin{equation}
    \partial_v \srfq (u, 0) + a(u) \partial_v \srfp (u, 0) = b(u) \partial_u \srfp (u, 0),
    \quad
    a(u)\ne  0,\
    u \in [0..1].
    \label{eq:g1}
\end{equation}
When
$b(u) := 0$ and $a(u) := const$, then we say the \splinepiece s join parameterically $C^1$.
Tensor-product and polar \sss s are internally $C^1$ (polar splines with a removable singularity at the pole, see \cite{Karc:2020:SPC}).
Otherwise the \splinepiece s join with geometric smoothness, short $G^1$.
For example, \cite{KKJP:2014:SMS}  uses 
a quadratic change of variables
to transition from the surrounding 
surface to one of the $n$  bi-3 patches of
the cap
and a linear one
between adjacent bi-3 patches of the cap.
 The result are explicit formulas that relate the input polyhedral net to the output BB-coefficients of the \sss, see
\cite{KKJP:2014:SMS}, \cite{sigT:}, \cite{Karc:2020:SPC}.
%
%
%
We note that
at global boundaries, as is typical for not-a-knot splines, the outermost layer recedes 
as in \figref{fig:bdy}.

\section{Computing with \sss s}
\label{sec:compute}

Bi-cubic \sss s have many applications
including industrial design,
visualization, animation, moment computation, re-approximation,
reconstruction, computation of partial differential equations on manifolds, etc..
For example, \cite{Nguyen:2016:CFE} uses a sub-class of \sss\ functions, for star-configurations only, to solve fourth order partial differential equations, and to compute geodesics on a free-form \sss\ surface via  the heat equation.

The appeal of \sss\ is that the applications are supported by the same representation as the geometry,
without additional meshing.
The irregular patterns can be viewed as 
both structurally necessary but also as a form of local adaptation.

First, in \secref{sec:test}, \ref{sec:impl} and \ref{sec:poisson-err}, we test the new elements 
by artificially introducing irregular polyhedral patterns into an
otherwise regular grid to be able to compare with a known exact solution.
The experiments indicate that the irregular patterns in \sss s do not noticeably increase the local error.
Second, in \secref{sec:freeform}, we compute second-order equations on a
free-form surface modeled by \sss s.

\subsection{Poisson Test}
\label{sec:test}
We evaluate the impact of inserting irregular mesh patterns on the error by solving a standard problem, 
Poisson's equation over the physical domain $\Omega = \srf{}(\unitdom)$, where $\unitdom :=
[-1..1]^2)$:
\begin{equation}
\text{find } \sol : \Omega \rightarrow \R :
\begin{cases}
	\Delta \sol  &= -f \hspace{3mm} \text{ in } \Omega, \\
	\hphantom{\Delta} \sol  &= 0 \hspace{7mm} \text{ on } \partial \Omega.
\end{cases}
\label{eq:poisson}
\end{equation}
The weak form  of Poisson's equation, projected into the $C^1$ space
of \sss\ basis functions $\Bps{}{j}$, is by Galerkin's approach,
\begin{align*}
   \int_\Omega \nabla
   \fh
   \cdot \nabla \Bps{}{i}(\idsrf{}) d_\Omega
   =
   \int_\Omega f\Bps{}{i}(\idsrf{}) d_\Omega,
   \quad
   \fh := \sum_j \bc_j \Bps{}{j} (\idsrf{}).
\end{align*}
The equation can be rewritten as the matrix equation
$K \bc = \mathbf{f}$ to be solved for the coefficient vector $\bc$ of $\fh$ where
\begin{align}
  K_{ij} &:= 
  \sum_{\pa}\int_{\unitdom} (\nabla \Bps{}{i})^\tr
      (J_{\pa}^{\tr} J_{\pa})^{-1}(\nabla \Bps{}{j})
       \df \unitdom,\quad
       J_\ga := \nabla _\st \srf{\ga},
      \\
      \mathbf{f}_i &:= \sum_{\pa} \int_\unitdom f \cdot \
      \Bps{i} \cdot \mathtt{J}\ \df\unitdom,
      \quad 
      \mathtt{J} := \sqrt{\deter (J_{\pa}^{\tr} J_{\pa})}.
  \label{eq:Galerkin_Lapl}
\end{align}
Here the sum is over all
pieces $\ga$ where $\Bps{}{i}$ has support.\\
\begin{wrapfigure}{r}{.5\linewidth} 
   \centering
   \begin{subfigure}[t]{0.60\linewidth}
      \centering
      \fcolorbox{lightgray}{white}{\includegraphics[trim=0 200 275 0, scale=0.22, clip=true]{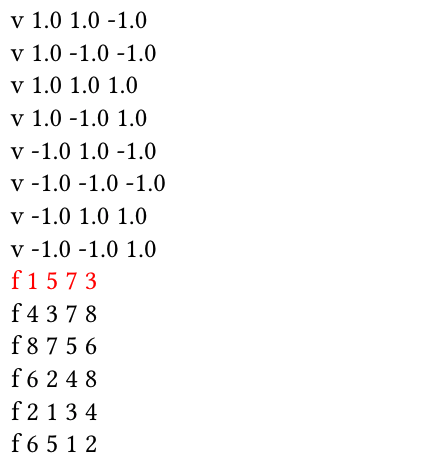}}
      \fcolorbox{lightgray}{white}
      {\includegraphics[trim=0 0 275 270, scale=0.22, clip=true]{fig/cube_obj_coef.png}}
      \caption{Input: cube.obj
      }
    \end{subfigure}
    \begin{subfigure}[t]{0.3
    \linewidth}
      \centering 
      \begin{overpic}[unit=1mm,scale=.25, width = \linewidth]{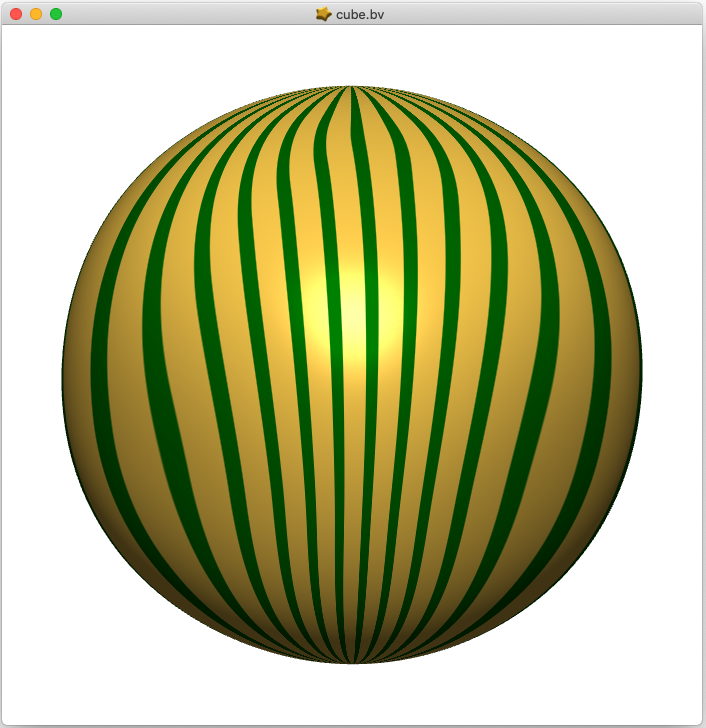}
      \put(-40,30) {$\Longrightarrow$}
      \end{overpic}
      \caption{surface} 
    \end{subfigure}
    \caption{\cite{toms:} The input .obj file and the output 
    \sss\ surface visualized  with \hl s in B\'ezierview \cite{Bview}.}
    \label{fig:io}
    \vskip-1cm
\end{wrapfigure}

\def\wid{0.22\linewidth}
\def\widS{0.22\linewidth}
\begin{figure}[]
    \centering
    \begin{subfigure}[t]{\widS}
        \includegraphics[width=\textwidth]{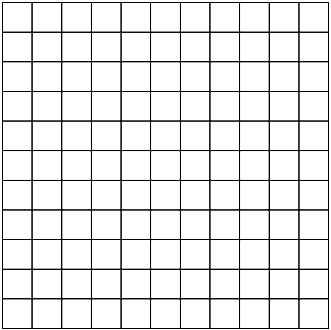}
    \end{subfigure} 
    \hspace{0.01\textwidth}
    \begin{subfigure}[t]{\wid}
        \includegraphics[width=\textwidth]{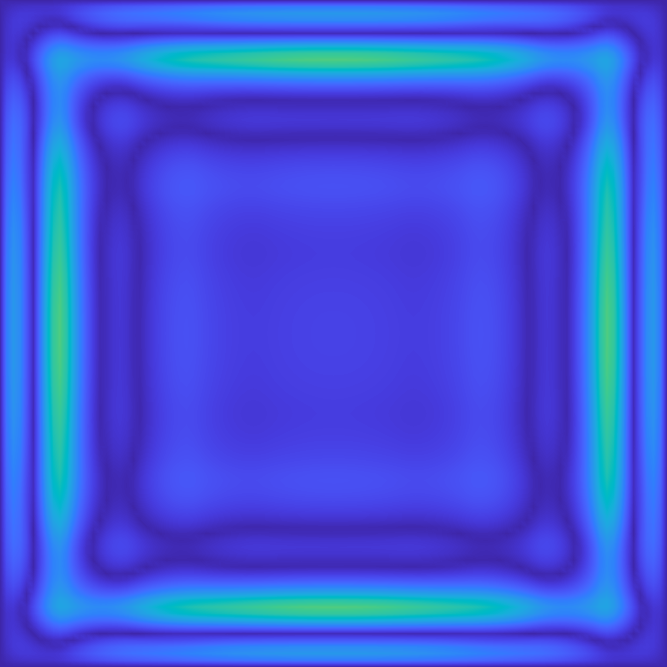}
        \subcaption{regular grid}
        \label{fig:regular-error}
    \end{subfigure}
    \skp
    \begin{subfigure}[t]{\widS}
        \includegraphics[width=\textwidth]{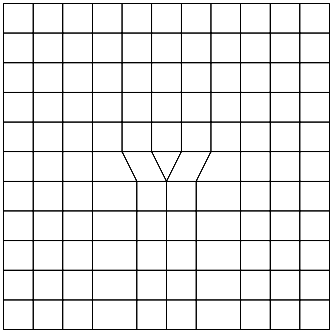}
    \end{subfigure} 
    \hspace{0.01\textwidth}
    \begin{subfigure}[t]{\wid}
        \includegraphics[width=\textwidth]{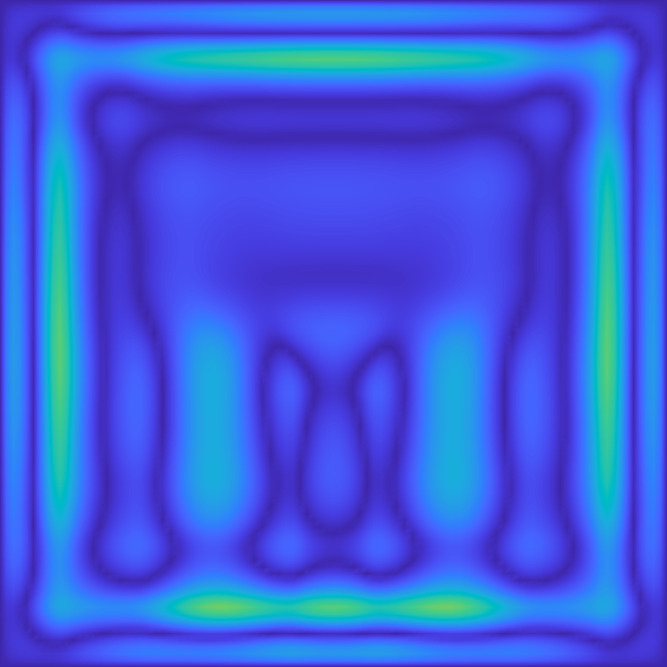}
        \subcaption{T0}
        \label{fig:t0-error}
    \end{subfigure} 
    \bigskip
    \\
    \begin{subfigure}[t]{\widS}
        \includegraphics[width=\textwidth]{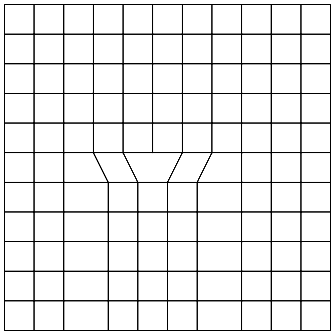}
    \end{subfigure} 
    \hspace{0.01\textwidth}
    \begin{subfigure}[t]{\wid}
        \includegraphics[width=\textwidth]{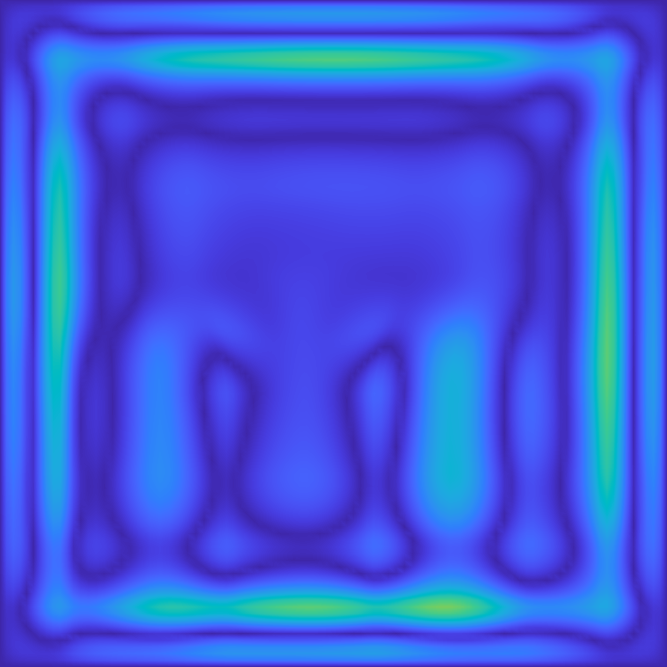}
        \subcaption{T1 subtract}
        \label{fig:t1-error-sub}
    \end{subfigure} 
    \skp
    \begin{subfigure}[t]{\widS}
        \includegraphics[width=\textwidth]{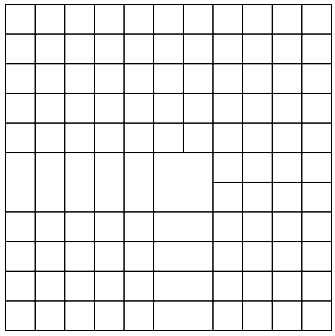}
    \end{subfigure} 
    \hspace{0.01\textwidth}
    \begin{subfigure}[t]{\wid}
        \includegraphics[width=\textwidth]{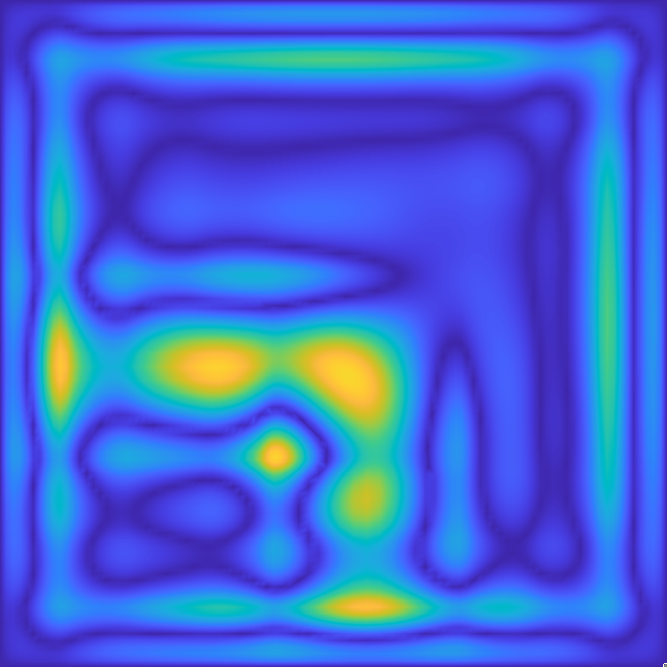}
        \subcaption{T2 subtract}
        \label{fig:t2-error-sub}
    \end{subfigure} 
    \bigskip
        \\
    \begin{subfigure}[t]{\widS}
        \includegraphics[width=\textwidth]{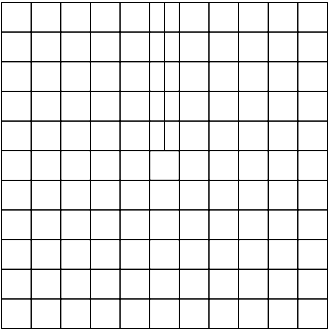}
    \end{subfigure} 
    \hspace{0.01\textwidth}
    \begin{subfigure}[t]{\wid}
        \includegraphics[width=\textwidth]{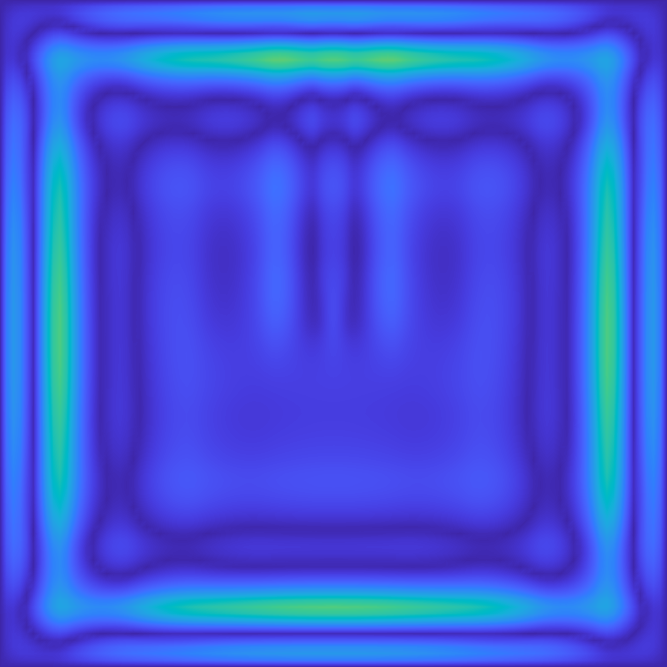}
        \subcaption{ T1 add}
        \label{fig:t1-error-add}
    \end{subfigure} 
    \skp
    \begin{subfigure}[t]{\widS}
        \includegraphics[width=\textwidth]{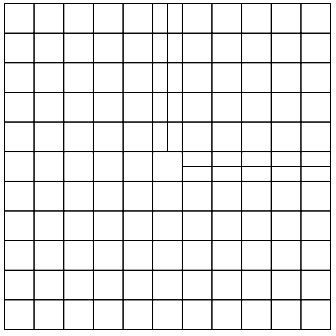}
    \end{subfigure}
    \hspace{0.01\textwidth}
    \begin{subfigure}[t]{\wid}
        \includegraphics[width=\textwidth]{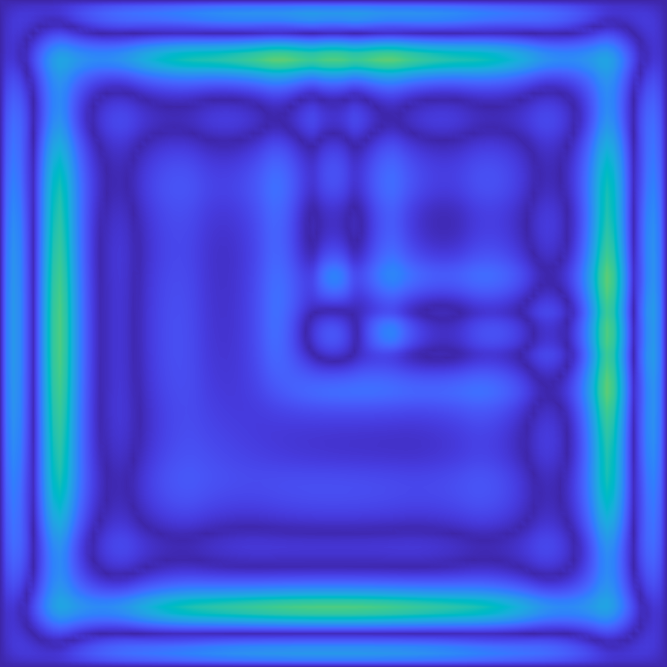}
        \subcaption{T2 add}
        \label{fig:t2-error-add}
    \end{subfigure}
    \bigskip
    \\
    \begin{subfigure}[t]{\widS}
        \includegraphics[width=\textwidth]{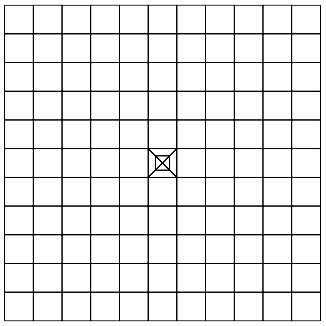}
    \end{subfigure} 
    \hspace{0.01\textwidth}
    \begin{subfigure}[t]{\wid}
        \includegraphics[width=\textwidth]{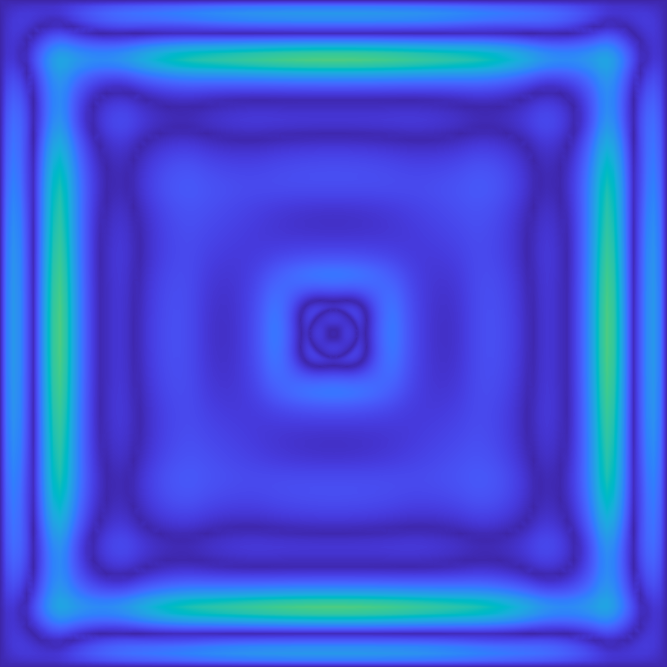}
        \subcaption{polar-net $n=4$}
        \label{fig:polar-net4}
    \end{subfigure} 
    \skp
    \begin{subfigure}[t]{\widS}
        \includegraphics[width=\textwidth]{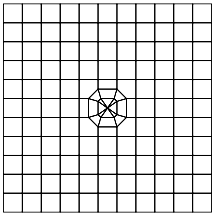}
    \end{subfigure}
    \hspace{0.01\textwidth}
    \begin{subfigure}[t]{\wid}
        \includegraphics[width=\textwidth]{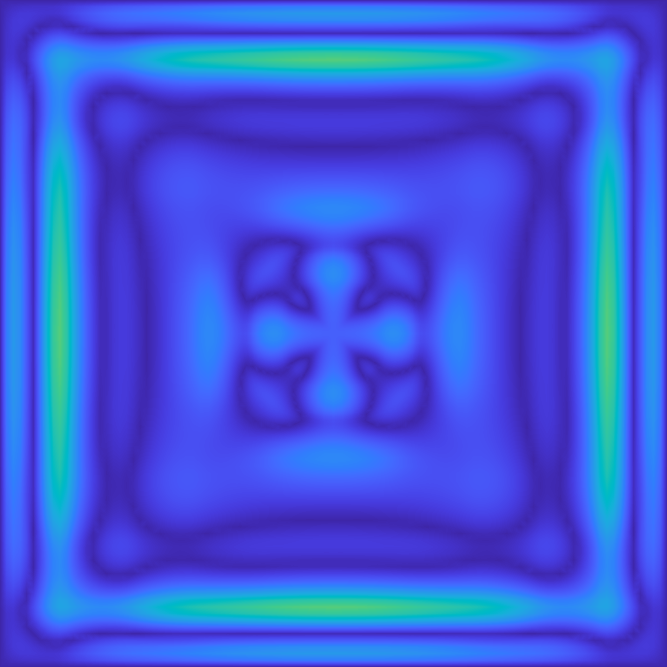}
        \subcaption{polar-net $n=8$}
        \label{fig:polar-net8}
    \end{subfigure}
       \\
    \bigskip
    \begin{subfigure}[t]{0.8\textwidth}
        \begin{overpic}[scale=.18,tics=10,width=\textwidth]{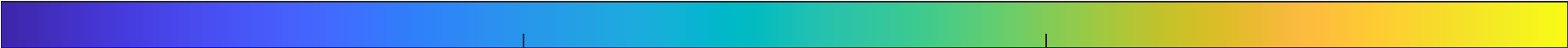}
            \put (0,-4) {0}
            \put (29,-4) {0.005}
            \put (63.5,-4) {0.01}
            \put (91.5,-4) {0.015}
        \end{overpic}
        \label{fig:error-bar}
    \end{subfigure}   
    \caption{
    Absolute error between approximate solution and exact solution for the various
    \sscn\ configurations. The exact solution $\sol = (u^2-1)(v^2-1)$
    varies between $0$ and $1$ on the domain
    $[-1..1]^2$.
    }
   \label{fig:sqPoissonErr}
\end{figure}
\subsection{Implementation Details}
\label{sec:impl}
The code \cite{toms:} collects the seven sub-net
configurations from the input faceted free-form  topological polyhedron
(in .obj format \cite{wavefront}, see \figref{fig:io}a) and generates the 
BB-coefficients of the pieces of a smooth spline
that can be visualized  by the online webGL viewer at \cite{Bview}, see \figref{fig:io}b.
Each $x$, $y$, $z$ coordinate is a \sss.
We intercept the output to compute the geometry terms 
$(J_{\pa}^{\tr} J_{\pa})^{-1}$ and $\mathtt{J}$
and access the connectivity information to determine where basis functions $\Bps{}{i}$
overlap in support to minimize the 
computational work in assembling $K$ and $\mathbf{f}$. Exact derivative 
polynomials are obtained by differencing 
BB-coefficients, e.g.\ $\bbc{{i+1},j} - \bbc{i,j}$ for differentiation with respect to the first variable.
For integration, all functions are evaluated at Gauss points,
shifted to the interval $[0..1]$ that is natural for the BB-form.  
Assembly is therefore no more costly than for  tensor-product splines.

We enforce the boundary constraints by 
modifying the basis functions
overlapping the boundary, setting to zero the 
BB-coefficients that determine   $\Bps{}{j}$ on the boundary  and beyond
(Recall that the surface boundary recedes by a half
from the boundary of the control polygon,
\figref{fig:bdy}.)
Here we interpret the outermost control points 
as belonging to bi-quadratic splines: they do, but have been degree-raised. 
Since for a bi-quadratic spline the only non-zero coefficient within the domain is on  the boundary, we can ignore their contribution.
We note that this approach may fail at corners for the $H^1$ norm since setting two partial derivatives to zero forces the derivatives at the corner to be zero.

\subsection{Error for Poisson's equation}
\label{sec:poisson-err}
To accurately measure the error, and so gauge
the impact of the irregular \sscn\ configurations, 
we choose a trivial geometric (physical surface) domain 
$x(u,v) = u$,
$y(u,v) = v$,
$z(u,v) = 0$  to cover $[-1..1]^2$. That is, 
we compute the graph of a function in
the knowledge that computations on more complex free-form surfaces are no more difficult
in the implemented framework (see \figref{fig:bottleheat}),
only harder to verify due to the lack of 
explicit solutions for free-form surfaces.
We choose 
$f := 2(u^2+v^2-2)$
so that we know 
 the exact solution of \eqref{eq:poisson} 
 to be
$\sol = (u^2-1)(v^2-1)$.

Starting with the tensor grid 
on $[-1..1]^2$ as a basic test, 
we solve and display the absolute difference between the \sss s and the exact solution of \eqref{eq:poisson} in 
\figref{fig:sqPoissonErr}.
Although the exact solution is bi-quadratic, enforcing the boundary constraints with fixed resolution results in an error along and near the boundary, of less than 1\%
see 
\figref{fig:regular-error}. Note the error bar
\figref{fig:sqPoissonErr}.
 \figref{fig:t1-error-sub} and
    \figref{fig:t2-error-sub} 
    are obtained by removing 
    grid-lines from 
            \figref{fig:regular-error}
whereas 
\figref{fig:t1-error-add} and
    \figref{fig:t2-error-add} 
    are obtained by adding
    grid-lines.
The location and disappearance of the
subtractive case
maximal 
errors
shows that the error is dominated by the reduction in the  
degrees of freedom when eliminating half a column 
and half a row of control points, and not due to distortion near the T-configuration.
 \figref{fig:t1-error-sub},
    \ref{fig:t2-error-sub},
\ref{fig:t1-error-add} and
\ref{fig:t2-error-add} 
remind of T-splines. 
However, T-splines \cite{sederberg2003t,kang2015new}
start with quadrilaterals and 
require a global parameterization
and then insert
partitions and T-junctions.
In general, the input free-form control nets cannot already include T-junctions
\cite[Fig 18]{karvciauskas2017t}.
\figref{fig:polar-net4} and \figref{fig:polar-net8} 
add minimal polar configurations with a decrease
in error near the pole and no noteworthy increase in error in the neighborhood.
    
We do not show
 convergence rates under refinement
because the \sss\ patterns do not have a natural refinement
and 
since the 2-norm 
error under refinement,
by splitting quadrilaterals into $2\times2$
or by applying Catmull-Clark subdivision to the 
\sscn,
is anyhow dominated by the regular mesh 
and 
the resulting spaces are not nested due to the different extent of the geometric continuity. To obtain nested refinement,
each quadrilateral bi-3 patch in BB-form can internally be partitioned and replaced by a $C^1$-connected spline complex that provides a local refinement hierarchy, the well-known T-spline space of `elements with hanging nodes', see e.g.\  \cite{kang2015new}.

\def\wid{0.25\linewidth}
\begin{figure}[]
    \centering
    \begin{subfigure}[t]{\wid}
        \includegraphics[width=\textwidth]{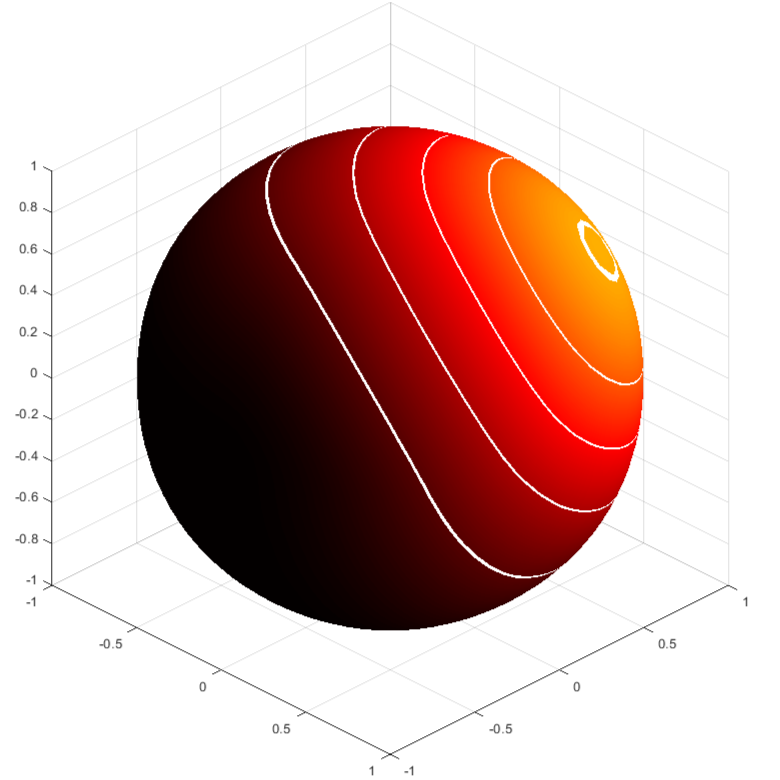}
        \subcaption{$t=0$}
        \label{fig:ch0}
    \end{subfigure} 
    \skp
    \begin{subfigure}[t]{\wid}
        \includegraphics[width=\textwidth]{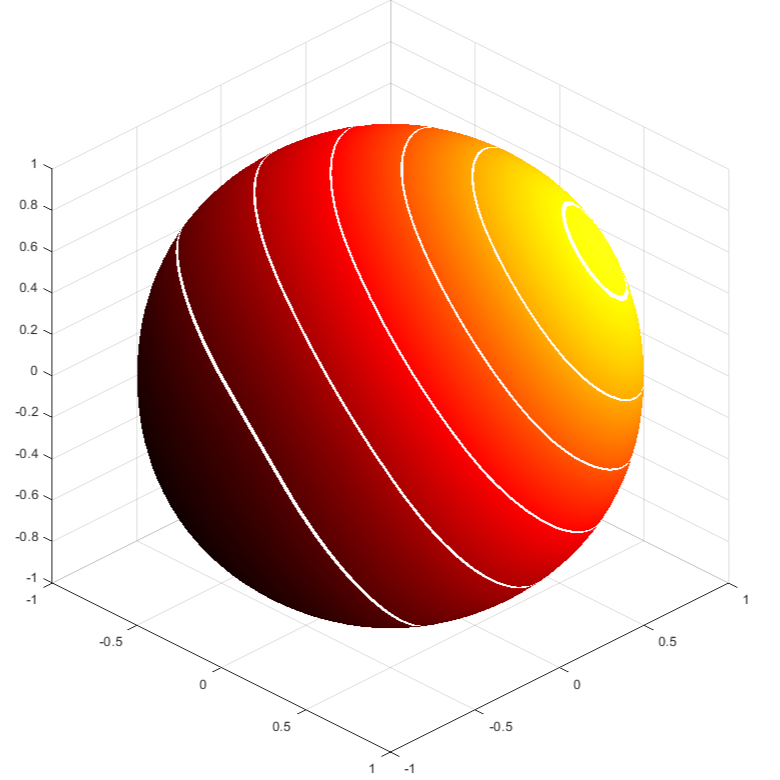}
        \subcaption{$t=0.25$}
        \label{fig:ch25}
    \end{subfigure}
    \skp
    \begin{subfigure}[t]{\wid}
        \includegraphics[width=\textwidth]{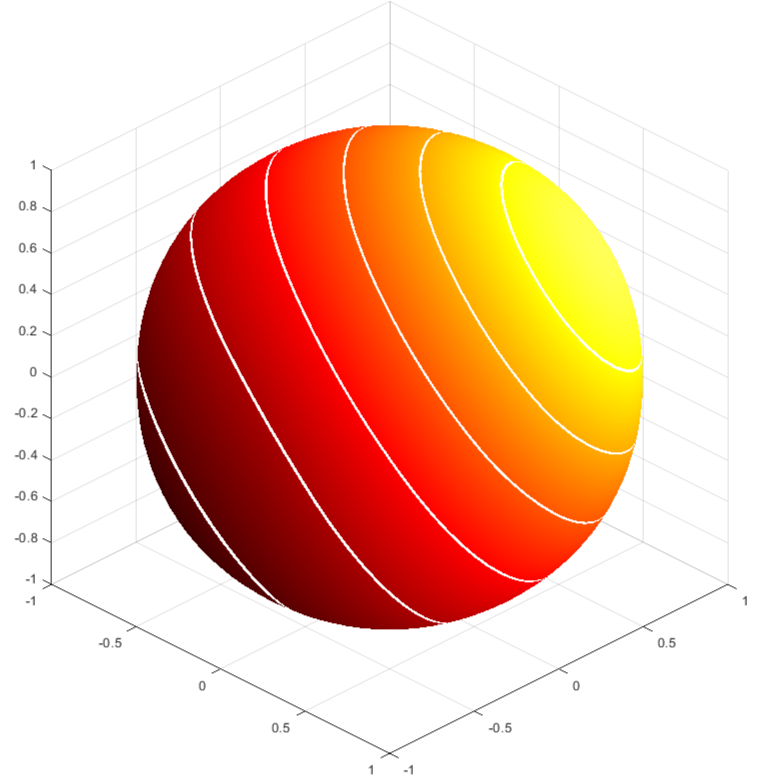}
        \subcaption{$t=0.5$}
        \label{fig:ch50}
    \end{subfigure}
    \skp
    \begin{subfigure}[t]{\wid}
        \includegraphics[width=\textwidth]{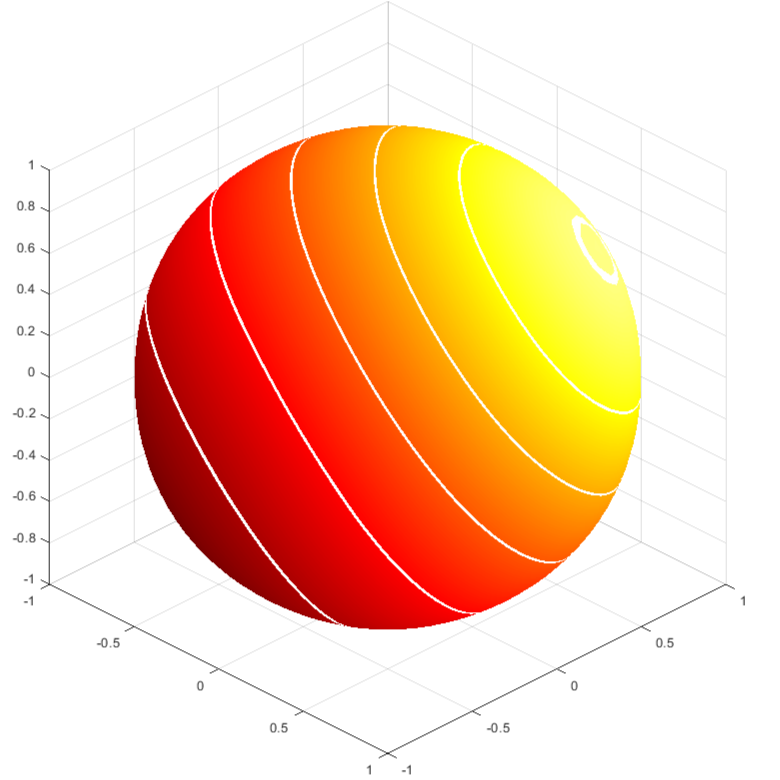}
        \subcaption{$t=0.75$}
        \label{fig:ch75}
    \end{subfigure}
    \skp
    \begin{subfigure}[t]{\wid}
        \includegraphics[width=\textwidth]{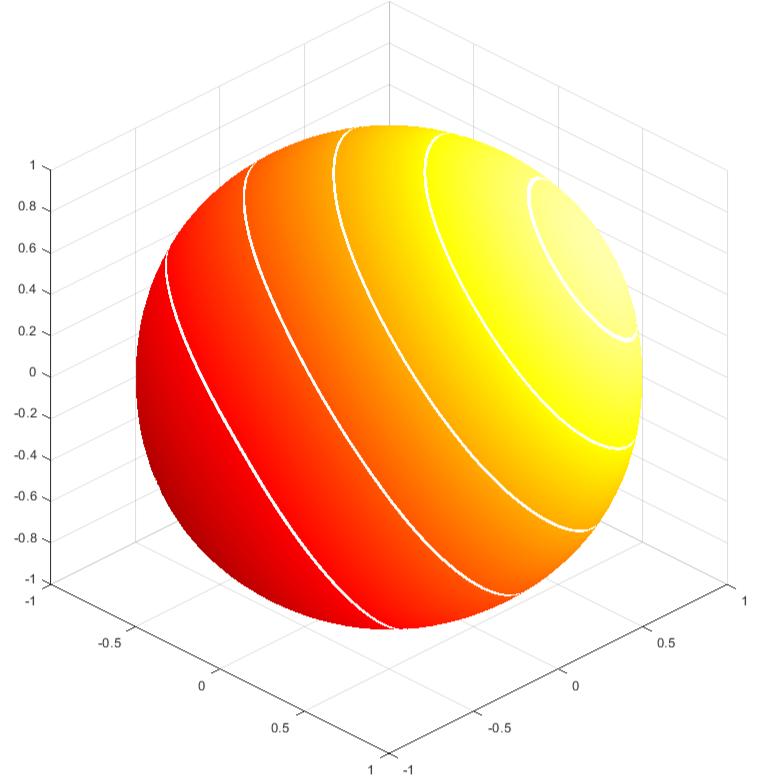}
        \subcaption{$t=1$}
        \label{fig:ch100}
    \end{subfigure}
    \bigskip
    \\
    \begin{subfigure}[t]{0.8\textwidth}
        \begin{overpic}[scale=.18,tics=10,width=\textwidth]{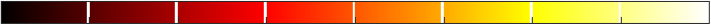}
            \put (0,-4) {0}
            \put (10.5,-4) {10}
            \put (23,-4) {20}
            \put (35.5,-4) {30}
            \put (48,-4) {40}
            \put (60.5,-4) {50}
            \put (73,-4) {60}
            \put (85.5,-4) {70}
            \put (96,-4) {80}
        \end{overpic}
    \end{subfigure}
   \bigskip
    \caption{
    Progression of geodesic fronts over time $t$ when solving the heat equation over the 
    polyhedral spline with input specified in \figref{fig:io}. The temperature 
    at the upper right
    vertex, the heat source,
    is kept constant to $80^\circ$.
    }
   \label{fig:time-series}
\end{figure}
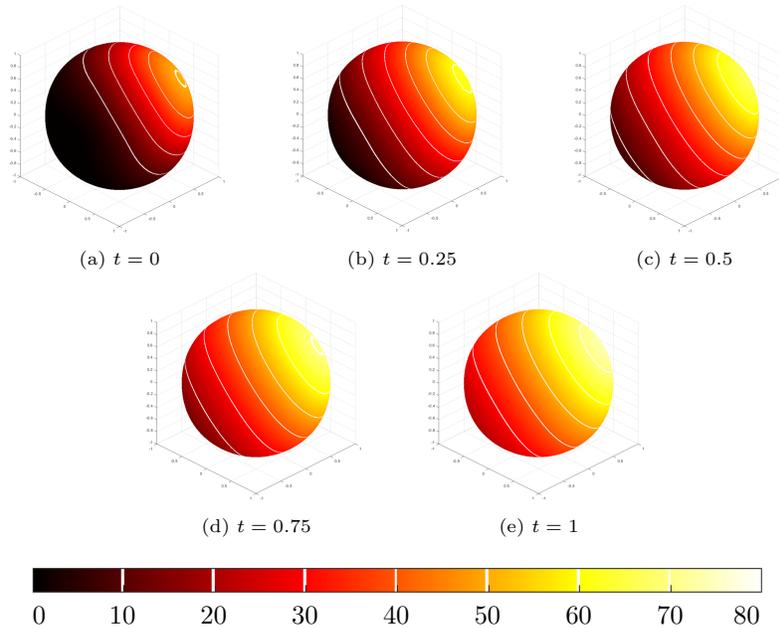
\subsection{Computing on free-form surfaces}
\label{sec:freeform}
Having verified correctness of the implementation
on the a square domain for all configurations,
we solve the heat equation on \sss\ free-form surfaces,
i.e.\ a second-order elliptic partial differential equation
with a solution evolving over time.
To certify correctness of the implementation, 
we start with the simple cube (8 nodes, 6 faces) of 
\figref{fig:io}.
We seed the temperature at one of the 8 vertices of the cube 
and observe the progress of the heat level curves over time.
This yields a progression of geodesic fronts,
as illustrated by the white lines in \figref{fig:time-series}. 
\def\wid{0.3\linewidth}
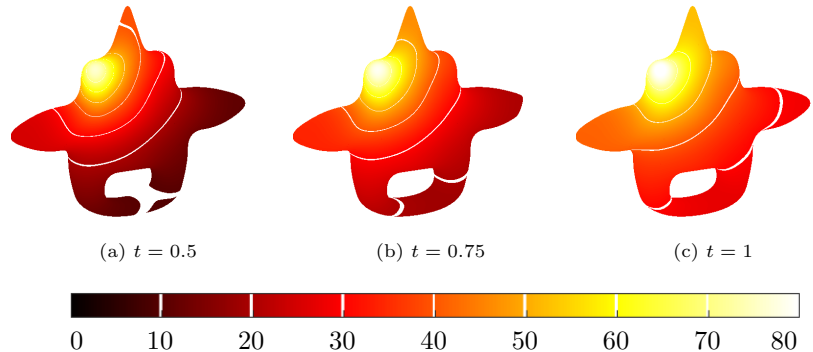
\begin{figure}[]
    \centering
    \begin{subfigure}[t]{\wid}
        \includegraphics[trim={80 80 0 0},clip,width=\textwidth]{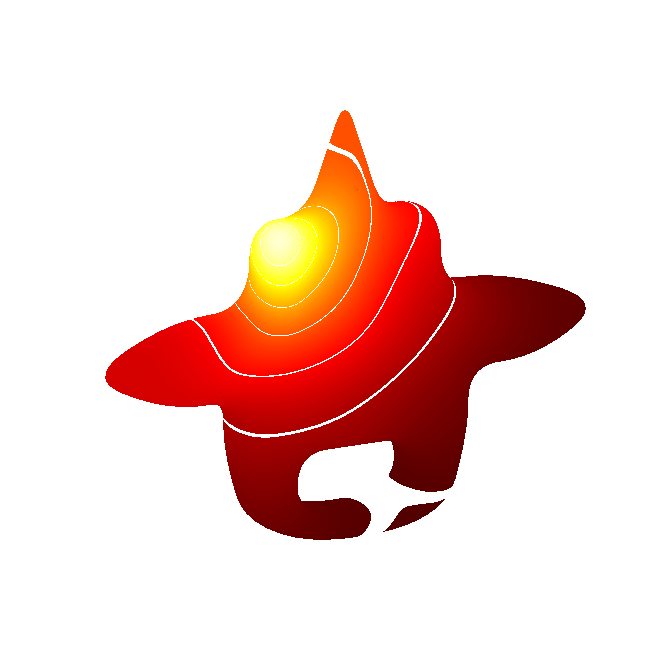}
        \subcaption{$t=0.5$}
        \label{fig:bh50}
    \end{subfigure}
    \begin{subfigure}[t]{\wid}
        \includegraphics[trim={80 80 0 0},clip,width=\textwidth]{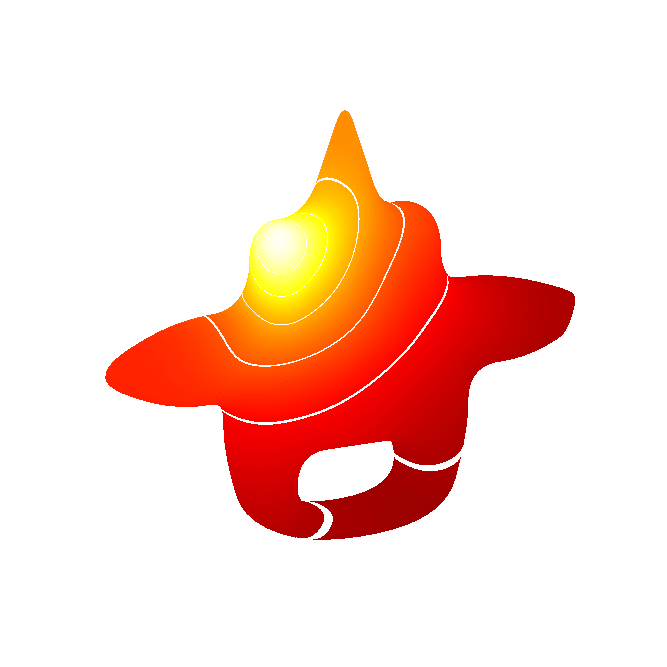}
        \subcaption{$t=0.75$}
        \label{fig:bh75}
    \end{subfigure}
    \begin{subfigure}[t]{\wid}
        \includegraphics[trim={80 80 0 0},clip, width=\textwidth]{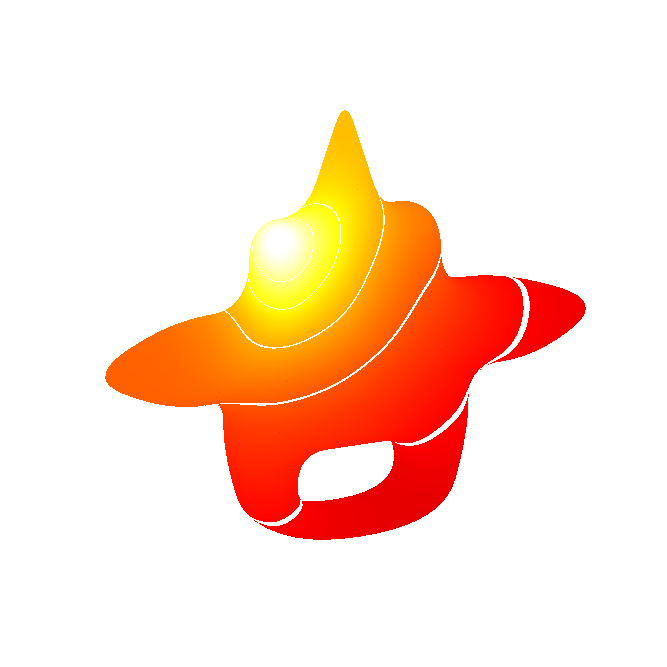}
        \subcaption{$t=1$}
        \label{fig:bh100}
    \end{subfigure}
    \bigskip
    \\
    \begin{subfigure}[t]{0.8\textwidth}
        \begin{overpic}[scale=.18,tics=10,width=\textwidth]{fig/colorbar_cube}
            \put (0,-4) {0}
            \put (10.5,-4) {10}
            \put (23,-4) {20}
            \put (35.5,-4) {30}
            \put (48,-4) {40}
            \put (60.5,-4) {50}
            \put (73,-4) {60}
            \put (85.5,-4) {70}
            \put (96,-4) {80}
        \end{overpic}
        \label{fig:error-bar}
    \end{subfigure}
    \caption{
    Progression of geodesic fronts for the 
   \sss\  with `bottle opener' input specified in \figref{fig:tightEnsemble} exposed to Florida sun
   at the upper left corner.
   (a) A low gradient causes a wide white level indicator curve at the lower right. The surface has no gaps.}
   \label{fig:bottleheat}
\end{figure}

Finally, to test a complex and  tight ensemble of polyhedral configurations, 
 we compute the heat progression 
on the`bottle opener' of \figref{fig:tightEnsemble}.
\figref{fig:bottleheat} shows the time series
before the opener becomes too hot to handle.

\section{Conclusion}
Analogous to tensor-product splines, but for more 
general control net configurations,
\sss s have been  shown to  model 
curved geometry (without trimming) and 
to represent higher-order functions on that geometry.
This setup has potential to be used for
engineering analysis on curved smooth objects
\emph{without additional meshing}.
As a proof-of-concept, we derived
time-evolving solutions of an elliptic partial
differential  equation. 
Since \sss s are differentiable,
also fourth-order equation Galerkin solvers 
make sense and can be implemented.

\bigskip
\noindent\textbf{Acknowledgement}
We thank William Gregory for support with C++ programming and Kyle Lo for generating input models
and explaining data structures of the \sss\ code.

\bibliographystyle{ACM-Reference-Format}
\bibliography{p}

\end{document}